\documentclass[11pt,epsfig]{article}
\usepackage{amsmath}
\usepackage{amsthm}
 \usepackage{amssymb}
\usepackage{epsfig}
 \usepackage{leftidx}
\usepackage{graphicx}
\usepackage{hyperref}
\usepackage{amsfonts}
\usepackage{graphicx, color, epstopdf}
\usepackage{verbatim}
\usepackage{CJK}
\usepackage{amsmath}

\newtheoremstyle{thmm}{1.5ex plus 1ex minus .2ex}{1.5ex plus 1ex minus .2ex}{\rmfamily}{}{\bfseries}{}{1em}{}
\theoremstyle{thmm}
\newtheorem{theorem}{Theorem}[section]
\newtheorem{lemma}{Lemma}[section]

\newtheorem{remark}{Remark}
\renewenvironment{proof}[1][Proof]{\noindent\textit{#1. } }{\hfill$\square$}
\def\a{\alpha}

\newcommand{\diff}{\triangledown_\tau}
\newcommand{\zd}{\,\mathrm{d}}

\newcommand{\absb}[1]{\big|#1\big|}

\newcommand{\bra}[1]{\left(#1\right)}
\newcommand{\brab}[1]{\big(#1\big)}
\newcommand{\braB}[1]{\Big(#1\Big)}
\newcommand{\brat}[1]{(#1)}
\newcommand{\kbra}[1]{\left[#1\right]}
\newcommand{\kbrab}[1]{\big[#1\big]}
\newcommand{\kbraB}[1]{\Big[#1\Big]}

\newcommand{\myinner}[1]{\left\langle#1\right\rangle}

\newcommand{\mynorm}[1]{\left\|#1\right\|}
\newcommand{\mynormb}[1]{\big\|#1\big\|}

\title{A fourth-order compact solver for fractional-in-time \\fourth-order diffusion equations}
\author{Jialing Zhong\thanks{Department of Mathematics, Nanjing University of Aeronautics and Astronautics,
211101, P. R. China. Jialing Zhong (jilzhong@163.com).}
\quad Hong-lin Liao\thanks{Corresponding author. Department of Mathematics,
Nanjing University of Aeronautics and Astronautics,
Nanjing 211106, P. R. China. Hong-lin Liao (liaohl@csrc.ac.cn,liaohl@nuaa.edu.cn)
is supported by a grant 1008-56SYAH18037
from NUAA Scientific Research Starting Fund of Introduced Talent.}
\quad Bingquan Ji\thanks{Department of Mathematics, Nanjing University of Aeronautics and Astronautics,
211101, P. R. China. Bingquan Ji (jibingquanm@163.com).}
\quad Luming Zhang\thanks{Department of Mathematics, Nanjing University of Aeronautics and Astronautics,
211101, P. R. China. Luming Zhang (zhanglm@nuaa.edu.cn)
is supported by the research grants No. 11571181 from National Natural Science Foundation of China.}
}
\begin{document}

\maketitle

\begin{abstract}
A fourth-order compact scheme is proposed for a fourth-order subdiffusion equation
with the first Dirichlet boundary conditions.  The fourth-order problem is firstly reduced into a couple of spatially
second-order system and we use an averaged operator to construct a fourth-order spatial approximation.
This averaged operator is compact since it involves only two grid points for the derivative boundary conditions.
The L1 formula on irregular mesh is considered for the Caputo fractional derivative,
so we can resolve the initial singularity of solution by putting more grid points near the initial time.
The stability and convergence are established by using three theoretical tools: a complementary discrete convolution kernel,
a discrete fractional Gr\"{o}nwall inequality and an error convolution structure.
Some numerical experiments are reported to demonstrate the accuracy and efficiency of our method.
\vspace{1em}

\noindent{\bf AMS subject classiffications.}\quad  65M06, 35B65
\end{abstract}

\section{Introduction}
\setcounter{equation}{0}
During the past several decades, fractional differential equations has become increasingly popular due to its wide applications in science and engineering \cite{Chen:2014,Hilfer:2000,MetzlerKlafter:2000},
including physics, chemistry, biochemistry and finance.
Despite analytic solutions of fractional differential equations may be found in
some special cases \cite{KilbasSrivastavaTrujillo:2006} by using Green function, Laplace and Fourier transforms,
most of practical problems can not be solved analytically.
Therefore, it is instructive to develop efficient numerical methods for time-fractional diffusion equations.
There are a lot of works contributed to the numerical solutions of subdiffusion and superdiffusion problems.
For examples, Sun and Wu \cite{SunWu:2006} constructed a fully discrete difference scheme by the method of order reduction for a diffusion-wave system,
and the corresponding  solvability, stability and  convergence were proved by the discrete energy method.
Based on the Gr\"{u}nwald-Letnikov discretization of Riemann-Liouville derivative,
Cui \cite{Cui:2009} developed and analyzed
a high-order compact finite difference scheme for solving one-dimensional fractional diffusion equation.
Gao and Sun \cite{GaoSun:2011} investigated a compact difference scheme
for the subdiffusion equation and proved the solvability, stability and convergence by the discrete energy method.

Apart from the second-order subdiffusion problem discussed in \cite{Cui:2009,GaoSun:2011,SunWu:2006} and the references therein,
the spatially fourth-order partial differential equations act an important role in modern science and engineering,
for instance, ice formation \cite{MyersCharpinChapman:2002,MyersCharpin:2004},
fluids on lungs \cite{HalpernJensenGrotberg:1998}
and the propagation of intense laser beams in a bulk medium with Kerr nonlinearity \cite{Karpman:1996}.
Also, there are a great amount of works on the fourth-order fractional
partial differential equations.
Agrawal \cite{Agrawal:2001} derived a general solution for a fourth-order fractional diffusion-wave equation
defined in a bounded space domain by using the finite sine transform technique and Laplace transform.
The solutions of
a generalized fourth-order fractional diffusion-wave equation was obtained in \cite{GolbabaiSayevand:2011} by
using the homotopy perturbation method.
Jafari et al. \cite{JafariDehghanSayevand:2008} showed that the Adomian decomposition method
is an useful analytical method for solving fourth-order fractional diffusion-wave equation.
One of typical fourth-order subdiffusion equations reads \cite{JiSunHao:2016},
\begin{eqnarray}\label{Problem-1}
\partial_t^{\alpha}u
+\frac{\partial^{4}u}{\partial{x}^{4}}
=qu+f(x,t)\quad\text{for $0<x<L$ and $0<t\leq T$,}
\end{eqnarray}
with an initial condition $u(x,0)=u_{0}(x)$.
Here, the reaction coefficient $q$ is a constant.
$\partial_t^{\alpha}u$ denotes the fractional Caputo derivative
of order $\alpha$ ($0<\alpha<1$) that is defined by \cite{Hilfer:2000},
\begin{eqnarray*}
(\partial_t^{\alpha}v)(t)
:=\int_{0}^{t} \omega_{1-\alpha} (t-s) \partial_s v(s)\zd{s}\quad\text{where}\quad
\omega_{\beta}(t):=\frac{t^{\beta-1}}{\Gamma(\beta)}\quad\text{for $t>0$}.
\end{eqnarray*}
In the literatures \cite{GuoLiDing:2014,HuZhang:2011,JiSunHao:2016,YaoWang:2018,ZhangPu:2017}, there are
several kinds of boundary conditions,  including
\begin{align*}
&\textrm{(BC0)}\quad
u(0,t)=b_{0l}(t),\quad u(L,t)=b_{0r}(t),\\
&\textrm{(BC1)}\quad
u_{x}(0,t)=b_{1l}(t),\quad u_{x}(L,t)=b_{1r}(t),\\
&\textrm{(BC2)}\quad
u_{xx}(0,t)=b_{2l}(t),\quad u_{xx}(L,t)=b_{2r}(t),\\
&\textrm{(BC3)}\quad
u_{xxx}(0,t)=b_{3l}(t),\quad u_{xxx}(L,t)=b_{3r}(t),
\end{align*}
where $0<t \leq T$.
Always, the combination (BC0)-(BC1) is called the first Dirichlet boundary conditions; (BC0) and (BC2) are called the second Dirichlet boundary conditions;
while (BC1)-(BC2), and the combination (BC1) and (BC3) are called Neumann boundary conditions.

In general, it is useful to introduce an auxiliary variable $v=\frac{\partial^{2}u}{\partial x^{2}}$ to transform the fourth-order subdiffusion equation
\eqref{Problem-1} into the following equivalent system
\begin{align}
\partial_t^{\alpha}u=&\,-\frac{\partial^{2}v}{\partial x^{2}}+qu+f(x,t)\quad\text{for $0<x<L$ and $0<t\leq T$,}\label{Problem-3}\\
v=&\,\frac{\partial^{2}u}{\partial x^{2}} \quad\text{for $0<x<L$ and $0\le t\leq T$,}\label{Problem-4}
\end{align}
subject to proper boundary conditions.

For treating the second Dirichlet boundary conditions (BC0) and (BC2),
one can use the classical compact operator $Aw_{i}:=\frac{1}{12}(w_{i-1}+10w_{i}+w_{i+1})$ for $1\leq i\leq M-1$,
because the boundary conditions for the variables $u$ and $v$ are Dirichlet-type
\cite{GuoLiDing:2014,HuZhang:2011,HuZhang:2012,ZhangPu:2017}.
A finite difference scheme for the fourth-order fractional diffusion-wave system with the second Dirichlet boundary conditions was proposed by Hu and Zhang \cite{HuZhang:2012}, and was proved to be uniquely solvable, stable and convergent in the $L^{\infty}$ norm by the discrete energy method.
They \cite{HuZhang:2011} also constructed a high-order compact difference scheme combining with the temporal extrapolation technique for the fourth-order fractional diffusion-wave system.
Guo et al. \cite{GuoLiDing:2014} derived two numerical schemes for a fourth-order subdiffusion equation. By using the Fourier method, they showed that the two finite difference schemes are unconditionally stable.
Zhang et al. \cite{ZhangPu:2017} proposed a compact scheme with a convergence order $O(\tau^{2}+h^{4})$ for the fourth-order subdiffusion equation with the second Dirichlet boundary conditions. Using the special properties of L2-1$_{\sigma}$ formula and the mathematical induction,
they proved the unconditional stability and convergence by discrete energy method.

However, the compact operators of the fourth-order derivative with other boundary conditions, such as first Dirichlet conditions
and Neumann boundary conditions, are quite different, especially at the boundary points.
Recently, Yao et al. \cite{YaoWang:2018} derived a compact difference scheme of order $O(\tau^{2}+h^{4})$
for fourth-order subdiffusion equations subject to Neumann boundary conditions (BC1) and (BC3).
The stability and convergence in the $L^{2}$ norm were established for the proposed scheme using
the following compact operator
\begin{align*}
Aw_{i}:=\left\{\begin{array}{cl}
\frac{5}{6}w_{0}+\frac{1}{6}w_{1}\,,& i=0,\vspace{2mm}\\
\frac{1}{12}(w_{i-1}+10w_{i}+w_{i+1})\,,&1\leq i\leq M-1,\vspace{2mm}\\
\frac{5}{6}w_{M}+\frac{1}{6}w_{M-1}\,,& i=M.
\end{array}\right.
\end{align*}
In addition, for the fourth-order subdiffusion equation \eqref{Problem-1} with the first Dirichlet boundary conditions
(BC0)-(BC1), Ji et al. \cite{JiSunHao:2016} developed a fourth-order scheme based on the following averaged operator
\begin{align*}
Aw_{i}:=\left\{\begin{array}{cl}
\frac{97}{180}w_{0}+\frac{19}{30}w_{1}-\frac{13}{60}w_{2}+\frac{2}{45}w_{3}\,,& i=0,\vspace{2mm}\\
\frac{1}{12}(w_{i-1}+10w_{i}+w_{i+1})\,,&1\leq i\leq M-1,\vspace{2mm}\\
\frac{97}{180}w_{M}+\frac{19}{30}w_{M-1}-\frac{13}{60}w_{M-2}+\frac{2}{45}w_{M-3}\,,& i=M.
\end{array}\right.
\end{align*}
The fully discrete scheme were constructed by combining the above operator for the spatial derivative
with the uniform L1 formula for the Caputo derivative. The difference scheme were proved to
be unconditionally stable and convergent in the $L^{2}$ norm by discrete energy method.
However, the spatial approximation is not compact since the averaged operator employs four grid points at the boundary points.

It is worth mentioning that the theoretical analysis and the corresponding convergence order in \cite{GuoLiDing:2014,HuZhang:2011,HuZhang:2012,JiSunHao:2016,YaoWang:2018,ZhangPu:2017} are always limited
because the solution of \eqref{Problem-1} is essentially nonsmooth near the initial time.
More seriously, as pointed out in \cite{RenLiaoZhangZhang}, the classical $H^{1}$ norm analysis and the $H^{2}$ norm analysis \cite{HuZhang:2011,HuZhang:2012} for the fourth-order problem always lead to
a loss of temporal accuracy when the solution is weakly singular near $t=0$.
In this article, we will construct a fourth-order compact difference scheme for
the fourth-order subdiffusion equation \eqref{Problem-1} with the  boundary conditions
(BC0)-(BC1), and establish sharp $L^2$ and $L^{\infty}$ norms error estimates under more realistic time regularity of solution.
The main contributions include:
\begin{itemize}
  \item[(1)] The new compact approximation of the boundary condition (BC1) involves only two grid points near the boundary, which is simpler than the approach \cite{JiSunHao:2016} using four grid points.
  \item[(2)] The initial singularity of solution is taken into account and resolved by employing nonuniform time steps.
      More interestingly, our method and the numerical analysis are available on general nonuniform meshes but not just some specific ones.
  \item[(3)] Sharp error estimates in the $L^2$ and maximum norms are obtained by applying an improved fractional Gr\"{o}nwall inequality
  and a convolution structure of consistency error.
\end{itemize}

We discretize the  time interval $[0,T]$ by
$0=t_{0}<t_{1}<t_{2}<\cdots <t_{N}=T$
with variable time-step sizes $\tau_{k}:=t_{k}-t_{k-1}$ for $1\leq k\leq N$.
Let the maximum step size $\tau:=\max_{1\leq k\leq N}\tau_{k}$,
and the adjoint step ratio $\rho_{k}:=\tau_{k}/\tau_{k+1}$ for $1\leq k\leq N-1$.
For the grid function $v:=\{v^{k}|0\leq k\leq N\}$,
let the difference operator $\diff v^k:=v^k-v^{k-1}$ for $k\geq1$.
For the numerical analysis of the fourth-order compact difference approximation in space, we
take $\Omega:=(0,L)$ and impose the following assumptions
\begin{align}
\mynormb{u(t)}_{H^{6}(\Omega)}\leq C_{u},\quad\mynormb{\partial_{t}u(t)}_{H^{6}(\Omega)}\leq C_{u}(1+t^{\sigma-1}),
\quad \mynormb{\partial_{tt}u(t)}_{H^{2}(\Omega)}\leq C_{u}(1+t^{\sigma-2})\label{regularityassumption-L2Norm}
\end{align}
for $0<t\leq T,$ where the parameter $\sigma\in(0,1)\cup(1,2)$ reflects the time regularity of solution.
To resolve the initial singularity of solution,
it is natural to put more mesh points near the initial time.
Specifically, the limitation of time steps is given as follows
\begin{itemize}
  \item[\textbf{AssG}] Let $\gamma \leq 1$ be a user-chosen parameter.
There is a constant $C_{\gamma}>0$, independent of $k$, such that $\tau_{k}\leq C_{\gamma}\tau \min\{1,t_{k}^{1-1/\gamma}\}$ for $1\leq k\leq N$ and $t_{k}\leq C_{\gamma}t_{k-1}$ for $2\leq k\leq N$.
\end{itemize}
Since $\tau_{1}=t_{1}$, \textbf{AssG} implies that $\tau_{1}=O(\tau^{\gamma})$. The parameter $\gamma$ controls the extent to which the grid points are concentrated near $t=0$. A practical example satisfying \textbf{AssG} is an initially graded grid $t_{k}=T(k/N)^{\gamma}$,  which has been discussed in \cite{LiaoLiZhang:2018,LiaoYanZhang:2018,StynesRiordanGracia:2016b}.

The remainder of this paper is organized as follows.
Some notations and auxiliary lemmas, are presented in next section.
Also, the fourth-order difference scheme and its numerical implementation are discussed in Section 2.
The stability and convergence of our method are established in Section 3.
Numerical examples are performed in Section 4 to demonstrate the accuracy and efficiency of the proposed scheme.
Throughout this article, any subscripted $C$, such as $C_{u}$ and $C_{v}$, denotes a generic positive constant,
not necessarily the same at different occurrences,
which may be dependent on the given data and the solution but independent of temporal and spatial mesh sizes.

\section{A fourth-order compact scheme}
\setcounter{equation}{0}
\subsection{Nonuniform L1 formula}
The well-known nonuniform L1 formula of Caputo derivative is denoted by
\begin{align}\label{L1formulaNonuniform}
D_{N}^{\alpha}v^{n}
:=\sum_{k=1}^n\int_{t_{k-1}}^{t_{k}}\omega_{1-\alpha}(t_n-s)\frac{\diff v^k}{\tau_k}\zd s
=\sum_{k=1}^na_{n-k}^{(n)}\diff v^k,
\end{align}
in which the discrete convolution coefficients $a_{n-k}^{(n)}$ are defined by
\begin{align}\label{Coefficient1-L1formulaNonuniform}
a_{n-k}^{(n)}:=\int_{t_{k-1}}^{t_{k}}\frac{\omega_{1-\alpha}(t_n-s)}{\tau_k}\zd s
=\frac{\omega_{2-\alpha}(t_n-t_{k-1})
-\omega_{2-\alpha}(t_n-t_{k})}{\tau_k}, \quad 1\leq k\leq n.
\end{align}
Furthermore, for fixed integer $n\geq 2$, the discrete L1 kernels of $a_{n-k}^{(n)}$ satisfy \cite{LiaoYanZhang:2018},
\begin{align}
a_{n-k-1}^{(n)}>\omega_{1-\alpha}(t_{n}-t_{k})>a_{n-k}^{(n)}, \quad 1\leq k \leq n-1.\label{Coefficient1-L1formulaNonuniform-monotone-1}
\end{align}

We recall the improved discrete fractional Gr\"{o}nwall inequality
in \cite{LiaoMcLeanZhang:2018,RenLiaoZhangZhang}, which is applicable for any nonuniform time meshes
and suitable for a variety of discrete fractional derivatives having a discrete form of
$\sum_{k=1}^na^{(n)}_{n-k}\triangledown_{\tau} v^k$.
The following fractional Gr\"{o}nwall lemma,
involving the well-known Mittag--Leffler function $E_\alpha(z):=\sum_{k=0}^\infty\frac{z^k}{\Gamma(1+k\alpha)},$
gathers previous (slightly simplified) results from Section 2 (including the main result, Theorems~2.3)
in \cite{RenLiaoZhangZhang}.
\begin{lemma}\label{FractGronwall}
Assume that the discrete convolution kernels $\{a_{n-k}^{(n)}\}_{k=1}^n$ satisfy the assumptions:
\begin{description}
\item[Ass1.] The discrete kernels are monotone, that is, $a_{k-2}^{(n)}\ge a_{k-1}^{(n)}>0$ for~$2\leq k\leq n\leq N$.
\item[Ass2.] There is a constant~$\pi_a>0$,
$a^{(n)}_{n-k}\ge\frac{1}{\pi_a\tau_k}\int_{t_{k-1}}^{t_k}\omega_{1-\alpha}(t_n-s)\zd{s}$
for $1\le k\le n\le N$.
\item[Ass3.] There is a constant~$\rho>0$ such that the local step ratio $\rho_k\le\rho$ for~$1\le k\le N-1$.
\end{description}
Define also a sequence of discrete complementary convolution kernels~$\{p_{n-j}^{(n)}\}_{j=1}^n$ by
\begin{align}\label{eq: discreteConvolutionKernel-RL}
p_{0}^{(n)}:=\frac{1}{a_0^{(n)}},\quad
p_{n-j}^{(n)}:=
\frac{1}{a_0^{(j)}}
\sum_{k=j+1}^{n}\brab{a_{k-j-1}^{(k)}-a_{k-j}^{(k)}}p_{n-k}^{(n)},
    \quad 1\leq j\leq n-1.
\end{align}
Then the discrete  complementary kernels $p^{(n)}_{n-j}\ge0$ are well-defined and fulfill
\begin{align}
&\sum_{j=k}^np^{(n)}_{n-j}a^{(j)}_{j-k}= 1
\quad\text{for $1\le k\le n\le N$.}\label{eq: P A}\\
&\sum_{j=1}^np^{(n)}_{n-j}\omega_{1+m\alpha-\alpha}(t_j)\leq \pi_{a}\omega_{1+m\alpha}(t_n)
\quad\text{for $m=0,1$ and $1\le n\le N$.}\label{eq: P bound}
\end{align}
Suppose that $\lambda$ is a non-negative constant independent of the time-steps
and the maximum step size $\tau\le1/\sqrt[\alpha]{2\Gamma(2-\alpha)\lambda\pi_{a}}.$
If the non-negative sequences $(v^k)_{k=0}^N$, $(\xi^{k})_{k=1}^{N}$ and $(\eta^{k})_{k=1}^{N}$ satisfy
\begin{equation}\label{eq: first Gronwall}
\sum_{k=1}^na^{(n)}_{n-k}\triangledown_{\tau}\brab{v^k}^2\le
    \lambda\brab{v^{n}}^2+v^{n}\xi^n+\bra{\eta^n}^2
    \quad\text{for $1\le n\le N$,}
\end{equation}
then it holds that, for $1\le n\le N$,
\begin{equation*}
v^n\le2E_\alpha\brab{2\max\{1,\rho\}\pi_a \lambda t_n^\alpha}
    \biggl(v^0+\max_{1\le k\le n}\sum_{j=1}^k p^{(k)}_{k-j}\xi^j+\sqrt{\pi_a\Gamma(1-\alpha)}\max_{1\le k\le n}\big\{t_k^{\alpha/2}\eta^k\big\}
    \biggr).
\end{equation*}
\end{lemma}

Note that, the definition \eqref{Coefficient1-L1formulaNonuniform} and the decreasing property
\eqref{Coefficient1-L1formulaNonuniform-monotone-1}
show that the discrete L1 kernels $a_{n-k}^{(n)}$ fulfill two assumptions
\textbf{Ass1}-\textbf{Ass2} in Lemma \ref{FractGronwall} with $\pi_{a}=1$. In what follows,
we will use the results
of Lemma \ref{FractGronwall}, including the complementary convolution kernels~$\{p_{n-j}^{(n)}\}_{j=1}^n$ defined by
\eqref{eq: discreteConvolutionKernel-RL}, without further declarations.

\subsection{Fourth-order compact approximation}
The discretization of space derivatives will be processed.
Take a positive integer $M$, 
let the spatial step size $h:=L/M$ and the discrete grid ${\Omega}_{h}:=\{x_{i}=ih\,|\,0\leq i \leq M\}$.
Define the function spaces on $\Omega_{h}$, $\mathcal{V}_{h}:=\big\{v\,|\,v=(v_{0},v_{1},\cdots v_{M})\big\}$ and
$
\mathring{\mathcal{V}}_{h}:=\big\{v\,|\,v\in\mathcal{V}_{h}\;\text{with}\; v_{0}=v_{M}=0\big\}.
$
For any function $v,w\in \mathcal{V}_{h}$,
define the difference operators
$\delta_{x}v_{i-\frac{1}{2}}:=(v_{i}-v_{i-1})/h$ and
\begin{eqnarray}\label{difference-operator-zero-one}
\delta_{x}^{2}v_{i}:=\left\{\begin{array}{cl}
\frac{2}{h}\delta_{x}v_{\frac{1}{2}}\,,& i=0,\vspace{2mm}\\
\frac{1}{h}(\delta_{x}v_{i+\frac{1}{2}}-\delta_{x}v_{i-\frac{1}{2}}),&1\leq i \leq M-1\,,\vspace{2mm}\\
-\frac{2}{h}\delta_{x}v_{M-\frac{1}{2}}\,,& i=M.
\end{array}\right.
\end{eqnarray}
As usual, define the discrete inner product
\begin{align*}
\myinner{v, w}:=\frac{h}{2}v_{0}w_{0}
+h\sum_{i=1}^{M-1}v_{i}w_{i}+\frac{h}{2}v_{M}w_{M},
\end{align*}
and the associated discrete $L^{2}$ norm $\|v\|:=\sqrt{\myinner{v, v}}$.
Also, we will use the discrete $H^1$ semi-norm $|u|_{1}:=\sqrt{\langle v, -\delta_{x}^2v\rangle}$,
the discrete $H^2$ semi-norm $\|\delta_{x}^{2}v\|:=\sqrt{\langle \delta_{x}^2v, \delta_{x}^2v\rangle}$
and the maximum norm  $\|v\|_{\infty}:=\max_{0\leq i \leq M}|v_{i}|$.
For any grid function $v\in \mathring{\mathcal{V}}_{h}$, by \cite[Lemma 2.5]{LiaoSun:2010} there exists a constant $C_{\Omega}>0$ such that
\begin{align}\label{EmbeddingInequalities}
\|v\|_{\infty}\leq C_{\Omega}\|\delta_{x}^{2}v\|.
\end{align}

In constructing the compact approximation, we employ the following averaged operator
\begin{align}\label{compact-operator-zero-one}
(\mathcal{A}v)_{i}:=\left\{\begin{array}{cl}
 \frac{2}{3}v_{0}+\frac{1}{3}v_{1}\,,& i=0,\vspace{2mm}\\
\frac{1}{12}(v_{i-1}+10v_{i}+v_{i+1})\,,&1\leq i\leq M-1,\vspace{2mm}\\
\frac{2}{3}v_{M}+\frac{1}{3}v_{M-1}\,,& i=M.
\end{array}\right.
\end{align}
Next two lemmas describe the local consistence errors of the averaged operator \eqref{compact-operator-zero-one}.
\begin{lemma}\cite[Lemma 4.1]{LiaoSun:2010}\label{lemma:discrete-space}
If $f$ is smooth and $\zeta(s):=5(1-s)^{3}-3(1-s)^{5}$, then
\begin{align*}
\frac{1}{12}\kbrab{f''(x_{i+1})+&10f''(x_{i})+f''(x_{i-1})}
=\frac{1}{h^{2}}\Big[f(x_{i+1})-2f(x_{i})+f(x_{i-1})\Big]\\
+&\,\frac{h^{4}}{360}\int_{0}^{1}\Big[f^{(6)}(x_{i}-sh)+f^{(6)}(x_{i}+sh)\Big]\zeta(s)\zd{s},\quad 1\leq i\leq M-1.
\end{align*}
\end{lemma}
\begin{lemma}\label{lemma:compact-operator-one-two-derivative-x0}
If $f$ is smooth and $\theta(s):=(1-s)^{3}\kbra{10-3(1-s)^{2}}$, then
\begin{align*}
\frac{2}{3}f''(x_{0})+\frac{1}{3}f''(x_{1})
=&\,\frac{2}{h}\kbraB{\frac{f(x_{1})-f(x_{0})}{h}-f'(x_{0})}
+\frac{h^2}{12}f^{(4)}(x_{0})\\
&\,+\frac{7h^3}{180}f^{(5)}(x_{0})
+\frac{h^4}{180}\int_{0}^{1}\theta(s)f^{(6)}(x_{0}+sh)\zd s,\\
\frac{2}{3}f''(x_{M})+\frac{1}{3}f''(x_{M-1})
=&\,\frac{2}{h}\kbraB{f'(x_{M})
-\frac{f(x_{M})-f(x_{M-1})}{h}}+\frac{h^2}{12}f^{(4)}(x_{M})\nonumber\\
&\,-\frac{7h^3}{180}f^{(5)}(x_{M})
+\frac{h^4}{180}\int_{0}^{1}\theta(s)f^{(6)}(x_{M}-sh)\zd s.
\end{align*}
\end{lemma}


\begin{proof}
The formula of Taylor expansion with integral remainder gives
\begin{align*}
f(x_{1})
=&\,f(x_{0})+hf'(x_{0})
+\frac{h^2}{2}f''(x_{0})
+\frac{h^3}{6}f'''(x_{0})
+\frac{h^4}{24}f^{(4)}(x_{0})\nonumber\\
&\,+\frac{h^5}{120}f^{(5)}(x_{0})
+\frac{h^6}{120}\int_{0}^{1}(1-s)^5f^{(6)}(x_{0}+sh)\zd s.
\end{align*}
A simple calculation shows that the above equality could be written as
\begin{align}\label{compact-operator-one-two-derivative-x0-proof-1}
\frac{2}{h}\kbraB{\frac{f(x_{1})-f(x_{0})}{h}&\,-f'(x_{0})}
=f''(x_{0})+\frac{h}{3}f'''(x_{0})+\frac{h^2}{12}f^{(4)}(x_{0})\nonumber\\
&\,+\frac{h^3}{60}f^{(5)}(x_{0})+\frac{h^4}{60}\int_{0}^{1}(1-s)^5f^{(6)}(x_{0}+sh)\zd s\,.
\end{align}
Again, using the Taylor expansion, we have
\begin{align}\label{compact-operator-zero-one-derivative-x0-proof-02}
f''(x_{1})
=&\,f''(x_{0})+hf'''(x_{0})+\frac{h^2}{2}f^{(4)}(x_{0})
+\frac{h^3}{6}f^{(5)}(x_{0})\nonumber\\
&\,+\frac{h^4}{6}\int_{0}^{1}(1-s)^3f^{(6)}(x_{0}+sh)\zd s.
\end{align}
Multiplying both sides of the above equality \eqref{compact-operator-zero-one-derivative-x0-proof-02} by $\frac{1}{3}$, and adding the factor $\frac{2}{3}f''(x_{0})$, we get
\begin{align}\label{compact-operator-zero-one-derivative-x0-proof-2}
\frac{2}{3}f''(x_{0})&\,+\frac{1}{3}f''(x_{1})=f''(x_{0})+\frac{h}{3}f'''(x_{0})+\frac{h^2}{6}f^{(4)}(x_{0})\nonumber\\
&\,+\frac{h^3}{18}f^{(5)}(x_{0})+\frac{h^4}{18}\int_{0}^{1}(1-s)^3f^{(6)}(x_{0}+sh)\zd s.
\end{align}
Subtracting  \eqref{compact-operator-one-two-derivative-x0-proof-1} from
\eqref{compact-operator-zero-one-derivative-x0-proof-2}
yields the first result. The second equality follows similarly and the proof is completed.
\end{proof}

\subsection{A fourth-order difference scheme}
Let the grid functions $u_{i}^{n}$ and $v_{i}^{n}$ be the discrete approximations of
exact
solutions $U_{i}^{n}:=u(x_{i},t_{n})$ and $V_{i}^{n}:=v(x_{i},t_{n})$,
respectively,
for $0\leq i\leq M$ and $0\leq n \leq N$.
From Lemma \ref{lemma:discrete-space}, it is easy to know that the compact approximations of \eqref{Problem-3}-\eqref{Problem-4} at the interior points read
\begin{align*}
\mathcal{A}D_{N}^{\alpha}u_{i}^{n}= &\,-\delta_{x}^{2}v_{i}^{n}+q\mathcal{A}u_{i}^{n}+\mathcal{A}f_{i}^{n},
\quad 1\leq i \leq M-1,1\leq n\leq N,\\
\mathcal{A}v_{i}= &\,\delta_{x}^{2}u_{i}^{n},\quad 1\leq i \leq M-1, 0\leq  n\leq N.
\end{align*}
For the left boundary point $x=0$, taking  the limit $x\rightarrow 0^{+}$ to the equation \eqref{Problem-3}, one has
\begin{align}\label{transtation1}
\partial_{xx}v(0,t)=(q-\partial_t^{\alpha})u(0,t)+f(0,t).
\end{align}
Moreover, we differentiate the governing equation \eqref{Problem-3} with respect to $x$ and find
\begin{align*}
\partial_{xxx}v(x,t)=(q-\partial_t^{\alpha})u_{x}(x,t)+f_{x}(x,t).
\end{align*}
Taking  the limit $x\rightarrow 0^{+}$ leads to
\begin{align}\label{transtation2}
\partial_{xxx}v(0,t)=(q-\partial_t^{\alpha})u_x(0,t)+f_{x}(0,t).
\end{align}
Combining the equations \eqref{transtation1}-\eqref{transtation2} with
Lemma \ref{lemma:compact-operator-one-two-derivative-x0},
we obtain a fourth-order approximation of the boundary condition (BC1) at the left boundary
\begin{align*}
\mathcal{A}v_{0}^{n}
=&\,\frac{2}{h}\brab{\delta_{x}u_{\frac{1}{2}}^{n}-b_{1l}(t_{n})}
+\frac{h^2}{12}\kbraB{(q-\partial_t^{\alpha})b_{0l}(t_{n})+f(0,t_{n})}\\
&\,+\frac{7h^3}{180}\kbraB{(q-\partial_t^{\alpha})b_{1l}(t_{n})+f_{x}(0,t_{n})
},\quad 0\leq n\leq N.
\end{align*}
The numerical approximation at the other boundary $x=L$ can be derived similarly.

In summary, we obtain the following compact scheme for the equations \eqref{Problem-3}-\eqref{Problem-4}
subject to boundary conditions (BC0)-(BC1),
\begin{align}
\mathcal{A}D_{N}^{\alpha}u_{i}^{n}=&\,-\delta_{x}^{2}v_{i}^{n}+q\mathcal{A}u_{i}^{n}+\mathcal{A}f_{i}^{n},
\quad 1\leq i \leq M-1,1\leq n\leq N,\label{scheme-zero-one-1}\\
\mathcal{A}v_{i}^{n}=&\,\delta_{x}^{2}u_{i}^{n},\quad 1\leq i \leq M-1, 0\leq n\leq N,\label{scheme-zero-one-2}\\
\mathcal{A}v_{0}^{n}=&\,\frac{2}{h}(\delta_{x}u_{\frac{1}{2}}^{n}
-b_{1l}(t_{n}))
+\frac{h^2}{12}\hat{b}_{0l}(t_{n})
+\frac{7h^3}{180}\hat{b}_{1l}(t_{n}),\quad 0\leq n\leq N,\label{scheme-zero-one-3}\\
\mathcal{A}v_{M}^{n}=&\,\frac{2}{h}(b_{1r}(t_{n})
-\delta_{x}u_{M-\frac{1}{2}}^{n})
-\frac{h^2}{12}\hat{b}_{0r}(t_{n})
+\frac{7h^3}{180}\hat{b}_{1r}(t_{n}),\quad 0\leq n\leq N,\label{scheme-zero-one-4}\\
u_{0}^{n}=&\,b_{0l}(t_{n}),\quad u_{M}^{n}=b_{0r}(t_{n}),\quad 1\leq n\leq N,\label{scheme-zero-one-5}
\end{align}
with $u_{i}^{0}=u_{0}(x_{i})$ for $0\leq i \leq M$,
where
\begin{align*}
&\hat{b}_{0l}(t_{n}):=(q-\partial_t^{\alpha})b_{0l}(t_{n})+f(0,t_{n})\,,\qquad
\hat{b}_{1l}(t_{n}):=(q-\partial_t^{\alpha})b_{1l}(t_{n})+f_{x}(0,t_{n})\,,\\
&\hat{b}_{0r}(t_{n}):=(q-\partial_t^{\alpha})b_{0r}(t_{n})+f(M,t_{n})\,,\quad
\hat{b}_{1r}(t_{n}):=(q-\partial_t^{\alpha})b_{1r}(t_{n})+f_{x}(M,t_{n})\,.
\end{align*}

\subsection{Numerical implementation}
The direct implementation of \eqref{scheme-zero-one-1}-\eqref{scheme-zero-one-4}
involves two independent variables $\{u_i^n,v_i^n\}$ and leads to a very large algebraic system of linear equations.
Here we eliminate the auxiliary grid function $\{v_{i}^{n}\}$ in \eqref{scheme-zero-one-1}-\eqref{scheme-zero-one-4} to obtain a self-contained difference system with respect to the original variables $\{u_i^n\}$.

Acting the difference operators $\mathcal{A}$ and $\delta_{x}^{2}$ on the equations \eqref{scheme-zero-one-1}-\eqref{scheme-zero-one-2} for $2\leq i \leq M-2$, respectively, and adding the resulting two equalities, we obtain the difference equation
\begin{align}
\mathcal{A}^{2}D_{N}^{\alpha}u_{i}^{n}+\delta_{x}^{4}u_{i}^{n}=q\mathcal{A}^{2}u_{i}^{n}
+\mathcal{A}^{2}f_{i}^{n},\quad 2\leq i \leq M-2,1\leq n\leq N.\label{scheme-equivalent-1}
\end{align}
According to the definition \eqref{compact-operator-zero-one} of $\mathcal{A}v_{0}^{n}$, one has
$$\mathcal{A}v_{0}^{n}\equiv h^{2}\brab{\frac{19}{36}\delta_{x}^{2}v_{1}^{n}
+\frac{1}{18}\delta_{x}^{2}v_{2}^{n}}+\frac{5}{3}\mathcal{A}v_{1}^{n}
-\frac{2}{3}\mathcal{A}v_{2}^{n}.$$
Using equations \eqref{scheme-zero-one-1} and \eqref{scheme-zero-one-2} for $i = 1, 2$,
we can rewrite the boundary scheme \eqref{scheme-zero-one-3} into
\begin{eqnarray}
&&\frac{19}{36}(\mathcal{A}D_{N}^{\alpha}u_{1}^{n}-q\mathcal{A}u_{1}^{n})
+\frac{1}{18}(\mathcal{A}D_{N}^{\alpha}u_{2}^{n}-q\mathcal{A}u_{2}^{n})\nonumber\\
&&\,\,+\frac{1}{h^2}\kbraB{\frac{2}{h}(\delta_{x}u_{\frac{1}{2}}^{n}-b_{1l}(t_{n}) )
+\frac{h^2}{12}\hat{b}_{0l}(t_{n})+\frac{7h^3}{180}\hat{b}_{1l}(t_{n}) -\frac{5}{3}\delta_{x}^{2}u_{1}^{n}+\frac{2}{3}\delta_{x}^{2}u_{2}^{n}}\,\nonumber\\
&&=\frac{19}{36}\mathcal{A}f_{1}^{n}+\frac{1}{18}\mathcal{A}f_{2}^{n}
,\quad 1\leq n\leq N.\label{scheme-equivalent-2}
\end{eqnarray}
Similarly, we have the following equation at the other boundary point
\begin{eqnarray}\label{scheme-equivalent-3}
&&\frac{19}{36}(\mathcal{A}D_{N}^{\alpha}u_{M-1}^{n}-q\mathcal{A}u_{M-1}^{n})
+\frac{1}{18}(\mathcal{A}D_{N}^{\alpha}u_{M-2}^{n}-q\mathcal{A}u_{M-2}^{n})
\nonumber\\
&&\,\, +\frac{1}{h^2}\kbraB{\frac{2}{h}(b_{1r}(t_{n})-\delta_{x}u_{M-\frac{1}{2}}^{n})
-\frac{h^2}{12}\hat{b}_{0r}(t_{n})+\frac{7h^3}{180}\hat{b}_{1r}(t_{n}) -\frac{5}{3}\delta_{x}^{2}u_{M-1}^{n}+\frac{2}{3}\delta_{x}^{2}u_{M-2}^{n}}\nonumber\\
&&=\frac{19}{36}\mathcal{A}f_{M-1}^{n}+\frac{1}{18}\mathcal{A}f_{M-2}^{n},\quad 1\leq n\leq N.
\end{eqnarray}

It is seen that the desired numerical solution $\{u_i^n\}$ can be computed by solving the linear difference equations \eqref{scheme-equivalent-1}
-\eqref{scheme-equivalent-3} together with the boundary values in \eqref{scheme-zero-one-5}.
Once the discrete solution $\{u_i^n\}$ is available, the auxiliary function $\{v_i^n\}$ can be obtained
by solving another self-contained algebraic system consisted of difference equations \eqref{scheme-zero-one-2}-\eqref{scheme-zero-one-4}.

\section{Stability and convergence}
\setcounter{equation}{0}

This section presents the numerical analysis of fourth-order difference
scheme \eqref{scheme-zero-one-1}-\eqref{scheme-zero-one-5}.
We introduce some preliminary lemmas, which are useful for the stability and
convergence analysis.

\begin{lemma}\cite[Lemma 5.2]{RenSun:2013}\label{lemma:compact-L2norm-equivalence}
For any grid function $u\in \mathring{\mathcal{V}}_{h}$, $\frac{1}{3}\|u\|^2\leq \|\mathcal{A}u\|^2\leq \|u\|^2.$
\end{lemma}

\begin{lemma}\label{lemma:compact-center-exchange}
For any grid functions $u, v\in \mathcal{V}_{h}$,
$\langle \delta_{x}^{2}v,\mathcal{A}u\rangle=\langle\delta_{x}^{2}u, \mathcal{A}v\rangle.$
\end{lemma}
\begin{proof} The definition \eqref{compact-operator-zero-one} of the operator $\mathcal{A}$ gives
\begin{align*}
(\mathcal{A}u)_{i}=\left\{\begin{array}{ll}
 u_{0}+\frac{h}{3}\delta_xu_{\frac12}\,,& i=0,\\
u_i+\frac{h^2}{12}\delta_x^2u_i\,,&1\leq i\leq M-1,\\
u_{M}-\frac{h}{3}\delta_xu_{M-\frac12}\,,& i=M.
\end{array}\right.
\end{align*}
Thus, by using the definition \eqref{difference-operator-zero-one} of $\delta_{x}^{2}$
and the discrete version of first Green formula
\begin{align*}
h\sum_{i=1}^{M-1}(\delta_{x}^{2}v_{i})u_i=-h\sum_{i=1}^{M}(\delta_{x}v_{i})(\delta_{x}u_{i})-u_0\delta_xv_{\frac12}+u_M\delta_xv_{M-\frac12},
\end{align*}
one derives that
\begin{align*}
\langle\delta_{x}^{2}v,\mathcal{A}u\rangle=&\,
\frac{h}{2}(\delta_{x}^{2}v_{0})(\mathcal{A}u_{0})
+h\sum_{i=1}^{M-1}(\delta_{x}^{2}v_{i})(\mathcal{A}u_{i})
+\frac{h}{2}(\delta_{x}^{2}v_{M})(\mathcal{A}u_{M})\nonumber\\
=&\,
\delta_{x}v_{\frac{1}{2}}\brab{u_{0}+\frac{h}{3}\delta_xu_{\frac12}}
+h\sum_{i=1}^{M-1}(\delta_{x}^{2}v_{i})\brab{u_i+\frac{h^2}{12}\delta_x^2u_i}
-\delta_{x}v_{M-\frac{1}{2}}\brab{u_{M}-\frac{h}{3}\delta_xu_{M-\frac12}}\nonumber\\
=&\,\frac{h}{3}(\delta_{x}v_{\frac{1}{2}})(\delta_{x}u_{\frac{1}{2}})
-h\sum_{i=1}^{M}(\delta_{x}v_{i})(\delta_{x}u_{i})
+\frac{h^3}{12}\sum_{i=1}^{M-1}(\delta_{x}^{2}v_{i})(\delta_{x}^{2}u_{i})
+\frac{h}{3}(\delta_{x}v_{M-\frac{1}{2}})(\delta_{x}u_{M-\frac{1}{2}})\nonumber\\
=&\,
\delta_{x}u_{\frac{1}{2}}\brab{v_{0}+\frac{h}{3}\delta_xv_{\frac12}}
+h\sum_{i=1}^{M-1}(\delta_{x}^{2}u_{i})\brab{v_i+\frac{h^2}{12}\delta_x^2v_i}
-\delta_{x}u_{M-\frac{1}{2}}\brab{v_{M}-\frac{h}{3}\delta_xv_{M-\frac12}}\nonumber\\
=&\,\frac{h}{2}(\delta_{x}^{2}u_{0})(\mathcal{A}v_{0})
+h\sum_{i=1}^{M-1}(\delta_{x}^{2}u_{i})(\mathcal{A}v_{i})
+\frac{h}{2}(\delta_{x}^{2}u_{M})(\mathcal{A}v_{M})=\langle\delta_{x}^{2}u,\mathcal{A}v\rangle.
\end{align*}
It completes the proof.
\end{proof}

\begin{lemma}\label{lemma:compact-fraction-exchange}
For any grid function $v^n\in \mathcal{V}_{h}$ for $0\le n\le N$,
\begin{align*}
D_{N}^{\alpha}\mathcal{A}v^{n}=\mathcal{A}D_{N}^{\alpha}v^{n}\quad\text{and} \quad
D_{N}^{\alpha}\delta_{x}^{2}v^{n}=\delta_{x}^{2}D_{N}^{\alpha}v^{n}.
\end{align*}
\end{lemma}

\begin{lemma}\label{lemma:time-equivalence}
For any function $v^n\in \mathcal{V}_{h}$ for $0\leq n\leq N$, $\langle D_{N}^{\alpha}v^{n},v^{n}\rangle \geq \frac12\sum_{k=1}^{n}a_{n-k}^{(n)}\nabla_{\tau}(\|v^{k}\|^{2}).$
\end{lemma}
\begin{proof}
The inequality can be derived from the proof of
\cite[Lemma 4.1]{LiaoMcLeanZhang:2018}.
\end{proof}

\subsection{Stability}
We present the stability of compact scheme \eqref{scheme-zero-one-1}-\eqref{scheme-zero-one-5} in the discrete $L^{2}$ and $L^{\infty}$ norms.
Consider the following perturbed system
\begin{align*}
\mathcal{A}D_{N}^{\alpha}\bar{u}_{i}^{n}=&\,-\delta_{x}^{2}\bar{v}_{i}^{n}
+q\mathcal{A}\bar{u}_{i}^{n}+\mathcal{A}f_{i}^{n}+\xi_{i}^{n}+\zeta_{i}^{n},
\quad 1\leq i \leq M-1,1\leq n\leq N,\\
\mathcal{A}\bar{v}_{i}^{n}=&\,\delta_{x}^{2}\bar{u}_{i}^{n}+\eta_i^n,\quad 1\leq i \leq M-1, 0\leq n\leq N,\\
\mathcal{A}\bar{v}_{0}^{n}=&\,\frac{2}{h}(\delta_{x}\bar{u}_{\frac{1}{2}}^{n}
-b_{1l}(t_{n}))
+\frac{h^2}{12}\hat{b}_{0l}(t_{n})
+\frac{7h^3}{180}\hat{b}_{1l}(t_{n})+\eta_0^n,\quad 0\leq n\leq N,\\
\mathcal{A}\bar{v}_{M}^{n}=&\,\frac{2}{h}(b_{1r}(t_{n})
-\delta_{x}\bar{u}_{M-\frac{1}{2}}^{n})
-\frac{h^2}{12}\hat{b}_{0r}(t_{n})
+\frac{7h^3}{180}\hat{b}_{1r}(t_{n})+\eta_M^n,\quad 0\leq n\leq N,\\
\bar{u}_{0}^{n}=&\,b_{0l}(t_{n}),\quad \bar{u}_{M}^{n}=b_{0r}(t_{n}),\quad 1\leq n\leq N,
\end{align*}
where $\xi_{i}^{n}$  and $\eta_{i}^{n}$ denote the exterior spatial forces,
while  $\zeta_{i}^{n}$ represents the exterior force in time direction.
 Let $\tilde{u}_i^n:=u_i^n-\bar{u}_i^n$ and $\tilde{v}_i^n:=v_i^n-\bar{v}_i^n$.
 We have the following perturbed equations of our numerical scheme \eqref{scheme-zero-one-1}-\eqref{scheme-zero-one-5},
\begin{align}
\mathcal{A}D_{N}^{\alpha}\tilde{u}_{i}^{n}
=&\,-\delta_{x}^{2}\tilde{v}_{i}^{n}+q\mathcal{A}\tilde{u}_{i}^{n}
+\xi_{i}^{n}+\zeta_{i}^{n}
,\quad 1\leq i \leq M-1 ,1\leq n\leq N,\label{scheme-stable-zero-one-1}\\
\mathcal{A}\tilde{v}_{i}^{n}
=&\,\delta_{x}^{2}\tilde{u}_{i}^{n}+\eta_{i}^{n},\quad 0\leq i \leq M,0\leq n\leq N,\label{scheme-stable-zero-one-2}
\end{align}
subject to the zero-valued boundary conditions $\tilde{u}_0^n=\tilde{u}_M^n=0$ for $1\leq n\leq N$.
Here and hereafter, $q_{+}:=\max\{q,0\}$ denotes the positive part of $q$.

\begin{theorem}\label{the:Stable-L2}
If the maximum time-step size $\tau\le1/\sqrt[\alpha]{4\Gamma(2-\alpha)q_{+}}$,
the discrete solution of perturbed equations \eqref{scheme-stable-zero-one-1}-\eqref{scheme-stable-zero-one-2} fulfills
\begin{align*}
\|\mathcal{A}\tilde{u}^n\|
\leq &\, 2E_{\alpha}\brab{4q_{+}\mbox{max}\{1,\rho\} t_n^{\alpha}}\braB{\|\mathcal{A}\tilde{u}^{0}\|
+2\max_{1\leq k\leq n}\sum_{j=1}^{k}p_{k-j}^{(k)}\|\zeta^{j}\|
\nonumber\\
&\,+2\Gamma(1-\alpha) \max_{1\leq k\leq n}\{t_{k}^{\alpha}\|\xi^{k}\|\}
+\sqrt{\Gamma(1-\alpha)} \max_{1\leq k\leq n}\{t_{k}^{\alpha/2 }\|\eta^{k}\|\}}\quad\text{for $1\leq n\leq N$. }
\end{align*}
So the compact scheme \eqref{scheme-zero-one-1}-\eqref{scheme-zero-one-5} is stable in the discrete $L^{2}$ norm.
\end{theorem}
\begin{proof}
Making the inner product of the equations \eqref{scheme-stable-zero-one-1}-\eqref{scheme-stable-zero-one-2} by
$2\mathcal{A}\tilde{u}^n$ and $2\mathcal{A}\tilde{v}^n$, respectively, and adding the two resulting equalities, one has
\begin{align}\label{scheme-stable-zero-one-proof-L2-1}
2\langle \mathcal{A}D_{N}^{\alpha}\tilde{u}^{n},\mathcal{A}\tilde{u}^{n}\rangle
+2\langle\mathcal{A}\tilde{v}^{n},\mathcal{A}\tilde{v}^{n}\rangle
=&\,-2\langle\delta_{x}^{2}\tilde{v}^{n},\mathcal{A}\tilde{u}^{n}\rangle
+2q\langle\mathcal{A}\tilde{u}^{n},\mathcal{A}\tilde{u}^{n}\rangle
+2\langle\xi^{n},\mathcal{A}\tilde{u}^{n}\rangle
\nonumber\\
&\,
+2\langle\zeta^{n},\mathcal{A}\tilde{u}^{n}\rangle
+2\langle\delta_{x}^{2}\tilde{u}^{n},\mathcal{A}\tilde{v}^{n}\rangle
+2\langle \eta^{n},\mathcal{A}\tilde{v}^{n}\rangle.
\end{align}
Lemma \ref{lemma:compact-center-exchange} shows that
$\langle\delta_{x}^{2}\tilde{v}^{n},\mathcal{A}\tilde{u}^{n}\rangle
=\langle\delta_{x}^{2}\tilde{u}^{n},\mathcal{A}\tilde{v}^{n}\rangle.
$
Obviously, the Cauchy-Schwarz and Young inequalities yield 
\begin{align*}
&2\langle\xi^{n},\mathcal{A}\tilde{u}^{n}\rangle\leq 2\|\mathcal{A}\tilde{u}^{n}\|\|\xi^{n}\|,\quad
2\langle\zeta^{n},\mathcal{A}\tilde{u}^{n}\rangle\leq 2\|\mathcal{A}\tilde{u}^{n}\|\|\zeta^{n}\|,\\
&2\langle \eta^{n},\mathcal{A}\tilde{v}^{n}\rangle \leq
\|\mathcal{A}\tilde{v}^{n}\|^{2}+\|\eta^{n}\|^{2}.
\end{align*}
Thus the equality \eqref{scheme-stable-zero-one-proof-L2-1} becomes
\begin{align}\label{scheme-stable-zero-one-proof-L2-4}
2\langle \mathcal{A}D_{N}^{\alpha}\tilde{u}^{n},\mathcal{A}\tilde{u}^{n}\rangle \leq &\,
2q\|\mathcal{A}\tilde{u}^{n}\|^{2}
+2\|\mathcal{A}\tilde{u}^{n}\|\bra{\|\xi^{n}\|+\|\zeta^{n}\|}
+\|\eta^{n}\|^{2}
\nonumber\\
\leq &\,
2q_{+}\|\mathcal{A}\tilde{u}^{n}\|^{2}
+2\|\mathcal{A}\tilde{u}^{n}\|\bra{\|\xi^{n}\|+\|\zeta^{n}\|}
+\|\eta^{n}\|^{2}.
\end{align}
Apply Lemma \ref{lemma:time-equivalence}
to the left hand side of \eqref{scheme-stable-zero-one-proof-L2-4}, it follows that
\begin{align*}
\sum_{k=1}^{n}a_{n-k}^{(n)}\nabla_{\tau}(\|\mathcal{A}\tilde{u}^{k}\|^{2}) \leq
2q_{+}\|\mathcal{A}\tilde{u}^{n}\|+2\|\mathcal{A}\tilde{u}^{n}\|\bra{\|\xi^{n}\|+\|\zeta^{n}\|}
+\|\eta^{n}\|^{2},
\end{align*}
which takes the form of \eqref{eq: first Gronwall} with the following substitutions
$$v^{k}:=\|\mathcal{A}\tilde{u}^{k}\|,\quad\lambda:=2q_{+},
\quad\xi^{n}:=2\bra{\|\xi^{n}\|+\|\zeta^{n}\|}\quad\text{and} \quad\eta^n:=\|\eta^{n}\|.$$
The discrete fractional Gr\"{o}nwall inequality in Lemma \ref{FractGronwall} gives the claimed inequality,
\begin{align*}
\|\mathcal{A}\tilde{u}^n\|
\leq &\, 2E_{\alpha}\brab{4q_{+}\mbox{max}\{1,\rho\} t_n^{\alpha}}\braB{\|\mathcal{A}\tilde{u}^{0}\|
+2\max_{1\leq k\leq n}\sum_{j=1}^{k}p_{k-j}^{(k)}\bra{\|\xi^{j}\|+\|\zeta^{j}\|}\nonumber\\
&\hspace{3cm}+\sqrt{\Gamma(1-\alpha)} \max_{1\leq k\leq n}\{t_{k}^{\alpha/2 }\|\eta^{k}\|\}}\quad \text{for $1\leq k\leq N.$}
\end{align*}
Thus the estimate \eqref{eq: P bound} of $m=0$ and Lemma \ref{lemma:compact-L2norm-equivalence} complete the proof.
\end{proof}

\begin{theorem}\label{the:Stable-Linfty}
If the maximum step size $\tau\le1/\sqrt[\alpha]{4\Gamma(2-\alpha)q_{+}}$,
then the numerical solution of of perturbed equations \eqref{scheme-stable-zero-one-1}-\eqref{scheme-stable-zero-one-2} satisfies
\begin{align}\label{scheme-stable-zero-one-conclusion-H2-1}
\|\delta_{x}^{2}u^{n}\|
\leq &\, 2E_{\alpha}(4q_{+}\mbox{max}\{1,\rho\} t_n^{\alpha})\braB{\|\delta_{x}^{2}u^{0}\|
+\|\eta^{0}\|+2\sum_{k=1}^{n}\|\diff \eta^{k}\|
+2\sqrt{3}\max_{1\leq k\leq n}\sum_{j=1}^{k}p_{k-j}^{(k)}\|\delta_x^2\zeta^{j}\|
\nonumber\\
&\,
+2q_{+}\Gamma(1-\alpha) \max_{1\leq k\leq n}\{t_{k}^{\alpha}\|\eta^{k}\|\}
+\sqrt{\Gamma(1-\alpha)} \max_{1\leq k\leq n}\{t_{k}^{\alpha/2 }\|\xi^{k}\|\}
}
+\|\eta^{n}\|
\end{align}
for $1\leq k\leq N$. So the numerical scheme \eqref{scheme-zero-one-1}-\eqref{scheme-zero-one-5} is stable in the $L^{\infty}$ norm.
\end{theorem}
\begin{proof}
Acting the difference operator $D_{N}^{\alpha}$ on \eqref{scheme-stable-zero-one-2} gives
\begin{align}
D_{N}^{\alpha}\mathcal{A}\tilde{v}_{i}^{n}
=&\,D_{N}^{\alpha}\delta_{x}^{2}\tilde{u}_{i}^{n}+D_{N}^{\alpha}\eta^{n},\quad 0\leq i \leq M,1\leq n\leq N.\label{scheme-stable-zero-one-3}
\end{align}
Taking the inner product of the equations \eqref{scheme-stable-zero-one-1} and \eqref{scheme-stable-zero-one-3} by $2\delta_x^2\tilde{v}^n$ and $2\mathcal{A}\tilde{v}^n$, respectively, and adding the two resulting equalities, we have
\begin{align}\label{scheme-stable-zero-one-proof-H2-0}
2\langle \mathcal{A}D_{N}^{\alpha}\tilde{u}^{n},\delta_x^2\tilde{v}^{n}\rangle
+2\langle D_{N}^{\alpha}\mathcal{A}\tilde{v}^{n},\mathcal{A}\tilde{v}^{n}\rangle
=-2\langle\delta_{x}^{2}\tilde{v}^{n},\delta_x^2\tilde{v}^{n}\rangle
+2q\langle\mathcal{A}\tilde{u}^{n},\delta_x^2\tilde{v}^{n}\rangle
+2\langle\xi^{n},\delta_x^2\tilde{v}^{n}\rangle\nonumber\\
+2\langle\zeta^{n},\delta_x^2\tilde{v}^{n}\rangle
+2\langle D_{N}^{\alpha}\delta_{x}^{2}\tilde{u}^{n},\mathcal{A}\tilde{v}^{n}\rangle
+2\langle D_{N}^{\alpha}\eta^{n},\mathcal{A}\tilde{v}^{n}\rangle.
\end{align}
Lemmas \ref{lemma:compact-center-exchange}-\ref{lemma:compact-fraction-exchange}
imply that
$$\langle \mathcal{A}D_{N}^{\alpha}\tilde{u}^{n},\delta_x^2\tilde{v}^{n}\rangle
=\langle D_{N}^{\alpha}\delta_{x}^{2}\tilde{u}^{n},
\mathcal{A}\tilde{v}^{n}\rangle
\quad \text{and} \quad
\langle\zeta^{n},\delta_x^2\tilde{v}^{n}\rangle
=\langle\delta_x^2\zeta^{n},\tilde{v}^{n}\rangle
.$$
Lemma \ref{lemma:compact-center-exchange} and the equation \eqref{scheme-stable-zero-one-2} yield
$$\langle \mathcal{A}\tilde{u}^{n}, \delta_{x}^{2}\tilde{v}^{n} \rangle
= \langle \delta_{x}^{2}\tilde{u}^{n}, \mathcal{A}\tilde{v}^{n} \rangle
= \|\mathcal{A}\tilde{v}^{n}\|^2
-\langle \eta^{n}, \mathcal{A}\tilde{v}^{n} \rangle.$$
Recalling the inequality $\frac{1}{3}\|u\|^{2}\leq \|\mathcal{A}u\|^{2}$ in Lemma \ref{lemma:compact-L2norm-equivalence}, we apply the Young and Cauchy-Schwarz inequalities to get
\begin{align*}
&2\langle\xi^{n},\delta_x^2\tilde{v}^{n}\rangle
\leq \|\delta_x^2\tilde{v}^{n}\|^{2}+\|\xi^{n}\|^{2}, \quad\quad\;
2\langle\delta_x^2\zeta^{n},\tilde{v}^{n}\rangle
\leq 2\sqrt{3}\|\mathcal{A}\tilde{v}^{n}\|\|\delta_x^2\zeta^{n}\|,
\nonumber\\
&2\langle D_{N}^{\alpha}\eta^{n},\mathcal{A}\tilde{v}^{n}\rangle
\leq 2\|\mathcal{A}\tilde{v}^{n}\|\|D_{N}^{\alpha}\eta^{n}\|,\quad
2\langle \eta^{n},\mathcal{A}\tilde{v}^{n}\rangle
\leq 2\|\mathcal{A}\tilde{v}^{n}\|\|\eta^{n}\|.
\end{align*}
Then, applying Lemma \ref{lemma:time-equivalence},
one obtains from \eqref{scheme-stable-zero-one-proof-H2-0} that
\begin{align*}
\sum_{k=1}^{n}a_{n-k}^{(n)}\nabla_{\tau}(\|\mathcal{A}\tilde{v}^{k}\|^{2})
\leq 2q_{+}\|\mathcal{A}\tilde{v}^{n}\|^{2}
+2\|\mathcal{A}\tilde{v}^{n}\|\brab{\sqrt{3}\|\delta_x^2\zeta^{n}\|
+\|D_{N}^{\alpha}\eta^{n}\|
+q_{+}\|\eta^{n}\|}
+\|\xi^{n}\|^{2},
\end{align*}
which has the form of \eqref{eq: first Gronwall} with the following substitutions $v^{k}:=\|\mathcal{A}\tilde{v}^{k}\|$,
$$\lambda:=2q_{+},\quad \xi^{n}:=2\sqrt{3}\|\delta_x^2\zeta^{n}\|
+2\|D_{N}^{\alpha}\eta^{n}\|+2q_{+}\|\eta^{n}\|\quad\text{and}\quad
\eta^{n}:=\|\xi^{n}\|.$$
The discrete fractional Gr\"{o}nwall inequality in Lemma \ref{FractGronwall} states that,
if the maximum time-step size
$\tau\le1/\sqrt[\alpha]{4\Gamma(2-\alpha)q_{+}}$, it holds that
\begin{align}\label{ieq:scheme-stable-zero-one-proof-H2-30}
\|\mathcal{A}\tilde{v}^{n}\|
\leq &\,2E_{\alpha}(4q_{+}\mbox{max}\{1,\rho\} t_n^{\alpha})\Big(\|\mathcal{A}\tilde{v}^{0}\|
+2\max_{1\leq k\leq n}\sum_{j=1}^{k}p_{k-j}^{(k)}\brab{\sqrt{3}\|\delta_x^2\zeta^{j}\|
+\|D_{N}^{\alpha}\eta^{j}\|
}
\nonumber\\
&\,+2q_{+}\max_{1\leq k\leq n}\sum_{j=1}^{k}p_{k-j}^{(k)}\|\eta^{j}\|
+\sqrt{\Gamma(1-\alpha)} \max_{1\leq k\leq n}\{t_{k}^{\alpha/2 }\|\xi^{k}\|\}\Big)\quad \text{for $1\leq n\leq N$}.
\end{align}
Applying the estimate \eqref{eq: P bound} of $m=0$, one has
\begin{align*}
\sum_{j=1}^{k}p_{k-j}^{(k)}\|\eta^{j}\|\le \Gamma(1-\a)\max_{1\leq j\leq k}\{t_{j}^{\a}\|\eta^{j}\|\}\quad\text{for $1\leq k\leq n$}.
\end{align*}
Using the L1 formula \eqref{L1formulaNonuniform},
we exchange the summation order to find that
\begin{align*}
\sum_{j=1}^{n}p_{n-j}^{(n)}\|D_{N}^{\alpha}\eta^{j}\| \leq \sum_{j=1}^{n}p_{n-j}^{(n)}\sum_{k=1}^{j}a_{j-k}^{(j)}\|\diff \eta^{k}\|
= \sum_{k=1}^{n}\|\diff \eta^{k}\|,
\end{align*}
where the identity \eqref{eq: P A} has been used in the equality.
Then one gets from \eqref{ieq:scheme-stable-zero-one-proof-H2-30} that
\begin{align*}
\|\mathcal{A}\tilde{v}^{n}\|
\leq &\, 2E_{\alpha}(4q_{+}\mbox{max}\{1,\rho\} t_n^{\alpha})\braB{\|\mathcal{A}v^{0}\|
+2\sqrt{3}\max_{1\leq k\leq n}\sum_{j=1}^{k}p_{k-j}^{(k)}\|\delta_x^2\zeta^{j}\|
+2\sum_{k=1}^{n}\|\diff \eta^{k}\|
\nonumber\\
&\,+2q_{+}\Gamma(1-\alpha) \max_{1\leq k\leq n}\{t_{k}^{\alpha}\|\eta^{k}\|\}
+\sqrt{\Gamma(1-\alpha)} \max_{1\leq k\leq n}\{t_{k}^{\alpha/2 }\|\xi^{k}\|\}
}
\end{align*}
for $1\leq k\leq N$. Then we employ the triangle inequality and the equation \eqref{scheme-stable-zero-one-2} to obtain
\begin{align*}
\|\delta_{x}^{2}\tilde{u}^{n}\| \leq \|\mathcal{A}\tilde{v}^{n}\|+\|\eta^{n}\|
\quad \text{and} \quad
\|\mathcal{A}\tilde{v}^{0}\| \leq \|\delta_{x}^{2}\tilde{u}^{0}\|+\|\eta^{0}\|.
\end{align*}
It yields the claimed estimate  \eqref{scheme-stable-zero-one-conclusion-H2-1}.
Then the embedding inequality \eqref{EmbeddingInequalities} implies that
the compact scheme \eqref{scheme-zero-one-1}-\eqref{scheme-zero-one-5} is stable in the discrete $L^{\infty}$ norm. It completes the proof.
\end{proof}

\subsection{Convergence}

Denote the local consistency error at time $t=t_n$ of the nonuniform L1 formula \eqref{L1formulaNonuniform} by
\[
\Upsilon^n[v]:=\partial_t^{\alpha}v(t_n)-D_{N}^{\alpha}v^{n},\quad n\geq1.
\]
Now we present the unconditional convergence of discrete solution in the discrete $L^{2}$ and $L^{\infty}$ norms.
It is to mention that, our convergence results are always valid on a general class of nonuniform meshes
(if the convergence order is not concerned),
because the error convolution structure of $\Upsilon^n[v]$
and the global consistency error $\sum^{n}_{j=1}p_{n-j}^{(n)}|\Upsilon^{n}[v]|$ in the next lemma are
valid without any priori information of time grids.
The detail proof can be found in Lemmas 3.1 and 3.3 (taking $\epsilon=0$) in \cite{LiaoYanZhang:2018}.

\begin{lemma}\cite[Lemmas 3.1 and 3.3]{LiaoYanZhang:2018}\label{lemma:local-consistency-error}
For $v\in C^2(0,T]$ with $\int_0^T t \,|v''(t)|\zd s < \infty$,
the local consistency error $\Upsilon^n[v]$ has the following error convolution structure
\begin{align*}
\absb{\Upsilon^{n}[v]}\leq a_{0}^{(n)}G^{n}+\sum_{k=1}^{n-1}(a_{n-k-1}^{(n)}-a_{n-k}^{(n)})G^{k}\quad \text{for $n\geq 1$},
\end{align*}
where $G^{k}$ is defined by
$$G^{k}:=2\int_{t_{k-1}}^{t_{k}}(t-t_{k-1})|v''(t)|\zd t\quad \text{for $1\le k\le n$.}$$
Suppose that there exists a constant $C_v>0$ such that $\absb{v''(t)}\leq C_v (1+t^{\sigma-2})$ for $0<t\leq T$,
where $\sigma\in(0,1)\cup(1,2)$ is a parameter. Then the global consistency error
\begin{align*}
\sum^{n}_{j=1}p_{n-j}^{(n)}\absb{\Upsilon^{j}[v]}\le 2\sum^{n}_{j=1}p_{n-j}^{(n)}a_{0}^{(j)}G^j\leq \frac{C_{v}}{\sigma}\tau_{1}^{\sigma}+\frac{C_{v}}{1-\alpha}\max_{2\leq k\leq n}t_{k}^{\alpha}t_{k-1}^{\sigma-2}\tau_{k}^{2-\alpha}\quad \text{for $n\geq 1$},
\end{align*}
where the discrete complementary convolution kernels $p_{n-j}^{(n)}$ are defined by \eqref{eq: discreteConvolutionKernel-RL}.
Specially, if the time mesh satisfies \textbf{AssG}, the global consistency error can be bounded by
\begin{align*}
\sum_{j=1}^{n}p_{n-j}^{(n)}\absb{\Upsilon^{j}[v]}\leq
\frac{C_{v}}{\sigma(1-\alpha)}\tau^{\mbox{min}\{2-\alpha,\, \gamma\sigma\}}\quad \text{for $n\geq 1$}.
\end{align*}
\end{lemma}

\begin{remark}
In the above global consistency error of L1 formula \eqref{L1formulaNonuniform},
the factor $\frac1{1-\alpha}$, tending to infinity as the fractional order $\alpha\rightarrow1$,
is mainly due to the application (taking $m=0$) of the rough estimate \eqref{eq: P bound} for
discrete complementary convolution kernels $p_{n-j}^{(n)}$.
It does not imply that the L1 formula \eqref{L1formulaNonuniform} can not employed when  $\alpha\rightarrow1$.
Actually, this factor disappears if we apply the case $m=1$ of \eqref{eq: P bound}
to evaluate the consistency error, although it would lead to a little lose of time accuracy.
\end{remark}

\begin{remark}[conjecture]
The discrete complementary convolution kernels $p_{n-j}^{(n)}$,
simulates the kernel of the
Riemann-Liouville integral
$(J_t^{\alpha}v)(t):=\int_{0}^{t} \omega_{\alpha} (t-s) v(s)\zd{s},$
see more details in \cite{LiaoLiZhang:2018,LiaoMcLeanZhang:2018} for the construction of $p_{n-j}^{(n)}$.
It is reasonable to conjecture that
\begin{align}\label{eq: P conjecture}
p_{n-j}^{(n)}\leq \pi_a\int_{t_{j-1}}^{t_{j}} \omega_{\alpha} (t_n-s)\zd{s}\quad \text{for $n\geq j\ge1$}.
\end{align}
because it directly makes the estimate \eqref{eq: P bound} available. Actually, we have
\begin{align*}
\sum_{j=1}^np^{(n)}_{n-j}&\, \omega_{1+m\alpha-\alpha}(t_j)\leq\pi_a\sum_{j=1}^n\int_{t_{j-1}}^{t_{j}} \omega_{\alpha} (t_n-s)\omega_{1+m\alpha-\alpha}(s)\zd{s}\\
=&\, \pi_a\int_{t_{0}}^{t_{n}} \omega_{\alpha} (t_n-s)\omega_{1+m\alpha-\alpha}(s)\zd{s}
=\pi_a\omega_{1+m\alpha}(t_n)
\quad\text{for $m=0,1$ and $n\ge1$.}
\end{align*}
In such case, one may derive a more sharp (pointwise) estimation of the global consistency error 
$\sum^{n}_{j=1}p_{n-j}^{(n)}\absb{\Upsilon^{j}[v]}$ of
nonuniform L1 formula \eqref{L1formulaNonuniform};
Nonetheless, up to now, we are not able to verify the estimate \eqref{eq: P conjecture}
from the definition \eqref{eq: discreteConvolutionKernel-RL} in mathematical manner.
\end{remark}

For the underlaying linear problem \eqref{Problem-1}, the essentially initial singularity can be resolved by using the graded-like time mesh
\textbf{AssG}. Let $e_i^n:=U_i^n-u_i^n$, $\epsilon_{i}^{n}:=V_i^n-v_i^n$ for $0\leq i\leq M$, $0\leq n\leq N$.
It is not difficult to find that the error functions $e_{i}^{n}\in \mathring{\mathcal{V}}_{h}$ and $\epsilon_{i}^{n}\in\mathcal{V}_{h}$ satisfy the following error system
\begin{align}
\mathcal{A}D_{N}^{\alpha}e_{i}^{n}=&\,
-\delta_{x}^{2}\epsilon_{i}^{n}
+q\mathcal{A}e_{i}^{n}
+(R_{su})_{i}^{n}+(\Upsilon^{n}[u])_i,\quad 1\leq i \leq M-1
, 1\leq n\leq N,\label{scheme-convegence-zero-one-1}\\
\mathcal{A}\epsilon_{i}^{n}=&\,\delta_{x}^{2}e_{i}^{n}+(R_{sv})_{i}^{n},\quad 0\leq i \leq M, 0\leq n\leq N,\label{scheme-convegence-zero-one-2}
\end{align}
where $R_{sv}$ and $R_{su}$ denote the truncation errors in space.

\begin{theorem}\label{the:Convergence-L2}
Suppose that the solution $u$ of problem \eqref{Problem-1} has the regularity property \eqref{regularityassumption-L2Norm} for the parameter $\sigma\in(0,1)\cup(1,2)$.
If the maximum time-step size $\tau\le1/\sqrt[\alpha]{4\Gamma(2-\alpha)q_{+}}$,
then the discrete solution of \eqref{scheme-zero-one-1}-\eqref{scheme-zero-one-5} is convergent with respect to the discrete $L^{2}$ norm, that is,
\begin{align*}
\|U^{n}-u^{n}\|\leq C_{u}\braB{\frac{\tau_{1}^{\sigma}}{\sigma}
+\frac{1}{1-\alpha}\max_{2\leq k\leq n}t_{k}^{\alpha}t_{k-1}^{\sigma-2}\tau_{k}^{2-\alpha}
+\frac{1}{1-\alpha}(t_{n}^{\alpha}
+t_{n}^{\alpha/2})h^4}\quad \text{for $1\leq n\leq N$.}
\end{align*}
Furthermore, if the mesh satisfies \textbf{AssG}, then
\begin{align}\label{scheme-convegence-zero-one-conclusion-L2-2}
\|U^{n}-u^{n}\|\leq \frac{C_{u}}{\sigma(1-\alpha)}\brab{\tau^{\min\{\gamma\sigma, 2-\alpha\}}+t_n^{\alpha}h^4}\quad \text{for $1\leq n\leq N$}.
\end{align}

\end{theorem}
\begin{proof}
By presenting a similar proof of Theorem \ref{the:Stable-L2}, we obtain that, for $1 \leq n \leq N$,
\begin{align}\label{scheme-convegence-zero-one-proof-L2-9}
\|\mathcal{A}e^n\| \leq&\, 2E_{\alpha}\brab{4q_{+}\mbox{max}\{1,\rho\} t_n^{\alpha}}
\braB{\|\mathcal{A}e^{0}\|+2\max_{1\leq k\leq n}\sum_{j=1}^{k}p_{k-j}^{(k)}\|\Upsilon^{j}[u]\|\nonumber\\
&\, +2\Gamma(1-\a)\max_{1\leq k\leq n}\{t_{k}^{\a}\|(R_{su})^{k}\|\}+\sqrt{\Gamma(1-\a)} \max_{1\leq k\leq n}\{t_{k}^{\a/2}\|(R_{sv})^{k}\|\}}.
\end{align}

We proceed to estimate the right-hand side of \eqref{scheme-convegence-zero-one-proof-L2-9}.
At first, $\|\mathcal{A}e^{0}\|=0$. Under the first regularity assumption in \eqref{regularityassumption-L2Norm}, one applies
Lemmas \ref{lemma:discrete-space} and \ref{lemma:compact-operator-one-two-derivative-x0}
to obtain the following spatial errors of fourth-order discretizations,
$$\|(R_{su})^{n}\|\leq C_{u}h^4\quad\text{for $1\leq n\leq N$} \quad\text{and}\quad\|(R_{sv})^{n}\|\leq C_{u}h^4 \quad \text{for $0\leq n\leq N$}.$$
By using Lemma \ref{lemma:local-consistency-error} combined with the third assumption in \eqref{regularityassumption-L2Norm},
the global consistency error
$$\sum_{k=1}^{n}p_{n-k}^{(n)}\|\Upsilon^{k}[u]\|\leq \frac{C_{u}}{\sigma}\tau_{1}^{\sigma}+\frac{C_{u}}{1-\alpha}\max_{2\leq k\leq n}t_{k}^{\alpha}t_{k-1}^{\sigma-2}\tau_{k}^{2-\alpha} \quad \text{for $1\leq n\leq N$}.$$
Thus, with the help of Lemma \ref{lemma:compact-L2norm-equivalence}, one obtains from \eqref{scheme-convegence-zero-one-proof-L2-9} that
\begin{align*}
\|U^{n}-u^{n}\|\leq C_{u}\braB{\frac{\tau_{1}^{\sigma}}{\sigma}+\frac{1}{1-\alpha}\max_{2\leq k\leq n}t_{k}^{\alpha}t_{k-1}^{\sigma-2}\tau_{k}^{2-\alpha}+\frac{1}{1-\alpha}(t_{n}^{\alpha}
+t_{n}^{\alpha/2})h^4}\quad \text{for $1\leq n\leq N$}.
\end{align*}
If the mesh satisfies \textbf{AssG}, it leads to the desired estimate \eqref{scheme-convegence-zero-one-conclusion-L2-2} and completes the proof.
\end{proof}
\begin{theorem}\label{the:Convergence-Linfty}
Assume that the solution $u$ of \eqref{Problem-1} fulfills the regularity assumption \eqref{regularityassumption-L2Norm} for the parameter $\sigma\in(0,1)\cup(1,2)$.
If the maximum time-step $\tau\le1/\sqrt[\alpha]{4\Gamma(2-\alpha)q_{+}}$,
the numerical solution of \eqref{scheme-zero-one-1}-\eqref{scheme-zero-one-5} is convergent
in the discrete $L^{\infty}$ norm, namely
\begin{align}\label{scheme-convegence-zero-one-conclusion-H2-1}
\|U^{n}-u^{n}\|_{\infty}\leq C_{u}\braB{\frac{\tau_{1}^{\sigma}}{\sigma}+\frac{1}{1-\alpha}\max_{2\leq k\leq n}t_{k}^{\alpha}t_{k-1}^{\sigma-2}\tau_{k}^{2-\alpha}
+\frac{1}{1-\alpha}\brab{1+t_{n}+t_{n}^{\alpha}+t_{n}^{\alpha/2}
+\frac{t_{n}^{\sigma}}{\sigma}
}h^4}
\end{align}
for $1 \leq n \leq N$. Specially, if the time mesh fulfills \textbf{AssG}, then
\begin{align}\label{scheme-convegence-zero-one-conclusion-H2-2}
\|U^{n}-u^{n}\|_{\infty}\leq \frac{C_{u}}{\sigma(1-\alpha)}\brab{\tau^{\min\{\gamma\sigma, 2-\alpha\}}+
\brat{t_{n}+t_{n}^{\alpha}+t_{n}^{\sigma}}h^{4}} \quad \text{for $1\leq n\leq N$}.
\end{align}
\end{theorem}
\begin{proof}
By presenting a similar proof of Theorem \ref{the:Stable-Linfty}, one obtains that, for $1\leq n\leq N$,
\begin{align}\label{scheme-convegence-zero-one-proof-H2-11}
\|\mathcal{A}\epsilon^{n}\|\leq &\, 2E_{\alpha}(4q_{+}\mbox{max}\{1,\rho\} t_n^{\alpha})\kbraB{\|\mathcal{A}\epsilon^{0}\|
+2\max_{1\leq k\leq n}\sum_{j=1}^{k}p_{k-j}^{(k)}
\brab{\sqrt{3}\|\delta_{x}^{2}\Upsilon^{j}[u]\|
+\|D_{N}^{\alpha}(R_{sv})^{j}\|
}
\nonumber\\
&\, +2q_{+}\Gamma(1-\alpha) \max_{1\leq k\leq n}\{t_{k}^{\alpha}\|(R_{sv})^{k}\|
+\sqrt{\Gamma(1-\alpha)} \max_{1\leq k\leq n}\{t_{k}^{\alpha/2}\|(R_{su})^{k}\|
\}}.
\end{align}

We proceed to estimate the right-hand side of \eqref{scheme-convegence-zero-one-proof-H2-11}.
With the help of the first assumption in \eqref{regularityassumption-L2Norm}, Lemmas \ref{lemma:discrete-space} and \ref{lemma:compact-operator-one-two-derivative-x0} imply that
$$\|(R_{su})^{n}\|\leq C_{u}h^4\quad\text{for $1\leq n\leq N$} \quad\text{and}\quad
\|(R_{sv})^n\|\leq C_{u}h^4\quad \text{for $0\leq n\leq N$}.$$
The error equation \eqref{scheme-convegence-zero-one-2} gives $\|\mathcal{A}\epsilon^{0}\|=\|(R_{sv})^0\|\le C_uh^4$.

Since the spatial error $(R_{sv})_i^n$ is defined uniformly
at the time $t=t_n$ (there is no temporal error in the equation \eqref{scheme-convegence-zero-one-2}),
we can define a time-continuous function $(R_{sv})_{i}(t)$ for $x_i\in\Omega_h$,
cf. section 4.3 in \cite{LiaoYanZhang:2018}, such that
$$(R_{sv})_{i}^n=(R_{sv})_{i}(t_n).$$
The second condition in \eqref{regularityassumption-L2Norm} implies $\mynorm{R_{sv}'(t)}\leq C_uh^4(1+t^{\sigma-1})$.
Hence, applying the L1 formula \eqref{L1formulaNonuniform} and the identity \eqref{eq: P A},
we exchange the summation order to find that
\begin{align*}
\sum_{j=1}^{n}p_{n-j}^{(n)}\|D_{N}^{\alpha}(R_{sv})^{j}\| \leq \sum_{j=1}^{n}p_{n-j}^{(n)}\sum_{k=1}^{j}a_{j-k}^{(j)}\|\diff (R_{sv})^{k}\|
= \sum_{k=1}^{n}\|\diff (R_{sv})^{k}\| \leq C_{u}h^4(t_{n}+\frac{t_{n}^{\sigma}}{\sigma}).
\end{align*}

Similarly, since the time consistency error $(\Upsilon^{n}[u])_i$ of L1 formula
can be defined uniformly with respect to the grid point $x_i\in\Omega_h$,
we define a space-continuous function $(\Upsilon^{n}[u])(x)$, cf. section 4.3 in \cite{LiaoYanZhang:2018},
such that $$(\Upsilon^{n}[u])_i=(\Upsilon^{n}[u])(x_i).$$
By using the Taylor expansion formula with integral remainder, we obtain
\begin{align*}
\delta_{x}^{2}(\Upsilon^{j}[u])(x_{i})=\int_{0}^{1}
\kbra{\partial_{xx}\Upsilon^{j}[u](x_{i}-sh)+\partial_{xx}\Upsilon^{j}[u](x_{i}+sh)}(1-s)
\zd s
\end{align*}
for $1 \leq i \leq M-1$ and $1 \leq j \leq n$. Lemma \ref{lemma:local-consistency-error}
with the third assumption in \eqref{regularityassumption-L2Norm} yields the global consistency error
\begin{align*}
\sum_{j=1}^{n}p_{n-j}^{(n)}\|\delta_{x}^{2}\Upsilon^{j}[u]\| \leq
C_{u}\braB{\frac{\tau_{1}^{\sigma}}{\sigma}+\frac{1}{1-\alpha}\max_{2\leq k\leq n}t_{k}^{\alpha}t_{k-1}^{\sigma-2}\tau_{k}^{2-\alpha}}.
\end{align*}

Therefore, collecting the above error estimates, one derives from
\eqref{scheme-convegence-zero-one-proof-H2-11} that
\begin{align*}
\|\mathcal{A}\epsilon^{n}\|\leq C_{u}\braB{\frac{\tau_{1}^{\sigma}}{\sigma}
+\frac{1}{1-\alpha}\max_{2\leq k\leq n}t_{k}^{\alpha}t_{k-1}^{\sigma-2}\tau_{k}^{2-\alpha}
+
\frac{1}{1-\alpha}\brab{t_{n}+t_{n}^{\alpha}+t_{n}^{\alpha/2}+\frac{t_{n}^{\sigma}}{\sigma}}h^4
}
\end{align*}
for $1\leq n\leq N$. Now we apply the triangle inequality and the error equation \eqref{scheme-convegence-zero-one-2} to get
\begin{align*}
\|\delta_x^2e^{n}\|\le \|\mathcal{A}\epsilon^{n}\|+\|(R_{sv})^n\|\le \|\mathcal{A}\epsilon^{n}\|+C_uh^4\quad \text{for $1\leq n\leq N$}
\end{align*}
Then the embedding inequality \eqref{EmbeddingInequalities} yields the error estimate \eqref{scheme-convegence-zero-one-conclusion-H2-1}.
When the time mesh satisfies \textbf{AssG}, the desired estimate \eqref{scheme-convegence-zero-one-conclusion-H2-2} follows from Lemma \ref{lemma:local-consistency-error}. The proof is completed.
\end{proof}


\section{Numerical experiments}
\setcounter{equation}{0}

We report numerical results to support the convergence theory.
The suggested compact difference scheme \eqref{scheme-zero-one-1}-\eqref{scheme-zero-one-5}
runs for solving the linear subdiffusion problem \eqref{Problem-1} on $\Omega=(0,1)$
over the time interval $(0,1]$.
We take $q=1$, $u_{0}(x)=0$, $b_{0l}(t)=0$, $b_{0r}(t)=0$, $b_{1l}(t)=\pi\omega_{1+\sigma}(t)$, $b_{1r}(t)=-\pi\omega_{1+\sigma}(t)$ and a source term
$f(x,t)=\kbra{\omega_{1+\sigma-\alpha}(t)+(\pi^4-1)\omega_{1+\sigma}(t)}\sin(\pi x).$
The exact solution of the subdiffusion problem \eqref{Problem-1} is $u=\omega_{1+\sigma}(t)\sin (\pi x)$.

Take a positive integer $N$, and consider a graded mesh $t_{k}=\bra{k/N}^{\gamma}T$ for $0\leq k \leq N$, where the grading parameter $\gamma\geq 1$ is chosen by the user.  The mesh is uniform if $\gamma=1$.
We employ an uniform spatial mesh with $M$ subintervals of length $h=L/M$. Furthermore, we measure the discrete maximum norm error
$e(M,N):=\max_{1\leq n \leq N}\|U(t_{n})- u^{n}\|_{\infty}$.
The spatial and temporal convergence rates are computed, respectively, by
$$\text{Order}(M)\approx\text{log}_{2}\bra{e(M,N)/e(2M,N)}\quad\text{and}\quad
\text{Order}(N)\approx\text{log}_{2}\bra{e(M,N)/e(M,2N)}.$$
The tests of spatial accuracy are reported in Tables 1--2, which confirm the fourth-order accuracy in space.
The temporal rate is examined in Tables 3--6 by four scenarios. The computational parameters are listed as follows,
\begin{itemize}
  \item Table 1: $N=10000$, $\sigma=1.3$ and $\gamma=2$ with fractional orders $\alpha=0.3$, $0.5$ and $0.7$.
  \item Table 2: $N=10000$, $\alpha=0.3$ and $\gamma=2$ with fractional orders $\sigma=1.3$, $1.5$ and $1.7$.
  \item Table 3: $M=100$, $\alpha=0.9$ and $\sigma=1.9$ with grid parameters $\gamma=1$, $1.5$ and $2$.

  \item Table 4: $M=600$, $\sigma=0.3$ and $\gamma=5$ with grid parameters $\alpha=0.3$, $0.5$, $0.7$.
  \item Table 5: $M=600$, $\alpha=0.5$ and $\sigma=0.3$ with grid parameters $\gamma=4$, $5$ and $6$.
  \item Table 6: $M=600$, $\alpha=0.4$ and $\sigma=0.3$ with grid parameters $\gamma=4$, $5$ and $6$.
\end{itemize}


\begin{table}[!ht]
\begin{center}
\tabcolsep 0pt {Table 1: \quad Numerical spatial accuracy for $N=10000$ and $\gamma=2$} \vspace*{0.5pt}
\def\temptablewidth{1.0\textwidth}
{\rule{\temptablewidth}{1pt}}
\begin{tabular*}{\temptablewidth}{@{\extracolsep{\fill}}ccccccc}
 $M$ &\multicolumn{2}{c}{$\alpha=0.3$,$\sigma=1.3$}&\multicolumn{2}{c}{$\alpha=0.5$,$\sigma=1.3$}
 &\multicolumn{2}{c}{$\alpha=0.7$,$\sigma=1.3$} \\
\cline{2-3}    \cline{4-5} \cline{6-7}& $e(M,N)$ & Order($M$) &$e(M,N)$ & Order($M$) &$e(M,N)$&Order($M$)\\
\hline	
8 &3.96e-04  &--  &3.96e-04&--  &3.96e-04&--  \\
16 &2.48e-05 &4.00&2.48e-05&4.00&2.48e-05&4.00\\
32 &1.55e-06 &3.99&1.55e-06&3.99&1.56e-06&3.99\\
64 &9.72e-08 &4.00&9.74e-08&4.00&9.93e-08&3.97 \\
\end{tabular*}
{\rule{\temptablewidth}{1pt}}
\end{center}\label{table:uniform1}
\end{table}

\begin{table}[!ht]
\begin{center}
\tabcolsep 0pt {Table 2: \quad Numerical spatial accuracy for $N=10000$ and $\gamma=2$} \vspace*{0.5pt}
\def\temptablewidth{1.0\textwidth}
{\rule{\temptablewidth}{1pt}}
\begin{tabular*}{\temptablewidth}{@{\extracolsep{\fill}}ccccccc}
 $M$ &\multicolumn{2}{c}{$\alpha=0.3$,$\sigma=1.3$}&\multicolumn{2}{c}{$\alpha=0.3$,$\sigma=1.5$}
 &\multicolumn{2}{c}{$\alpha=0.3$,$\sigma=1.7$} \\
\cline{2-3}    \cline{4-5} \cline{6-7}& $e(M,N)$ & Order($M$) &$e(M,N)$ & Order($M$) &$e(M,N)$&Order($M$)\\
\hline		
8 &3.96e-04  &--  &3.47e-04&--  &2.99e-04&--  \\
16 &2.48e-05 &4.00&2.17e-05&4.00&1.87e-05&4.00\\
32 &1.55e-06 &3.99&1.36e-06&3.99&1.17e-06&3.99\\
64 &9.72e-08 &4.00&8.53e-08&4.00&7.34e-08&4.00 \\
\end{tabular*}
{\rule{\temptablewidth}{1pt}}
\end{center}\label{table:uniform1}
\end{table}

\begin{table}[!ht]
\begin{center}
\tabcolsep 0pt {Table 3: \quad Numerical temporal accuracy for
$M=100,\alpha=0.9$ and $\sigma=1.9$ } \vspace*{0.5pt}
\def\temptablewidth{1.0\textwidth}
{\rule{\temptablewidth}{1pt}}
\begin{tabular*}{\temptablewidth}{@{\extracolsep{\fill}}ccccccc}
 $N$ &\multicolumn{2}{c}{$\gamma=1$}
 &\multicolumn{2}{c}{$\gamma=1.5$}&\multicolumn{2}{c}{$\gamma=2$} \\
\cline{2-3}    \cline{4-5} \cline{6-7}& $e(M,N)$ & Order($N$)
&$e(M,N)$ & Order($N$) &$e(M,N)$&Order($N$) \\
\hline
128 &6.56e-06 &--&6.21e-06&--&7.37e-06&--\\
256 &3.16e-06 &1.05&2.90e-06&1.10&3.41e-06&1.11 \\
512&1.52e-06 &1.06&1.36e-06&1.10&1.59e-06&1.10 \\
1024&7.32e-07 &1.07& 6.36e-07&1.09&7.42e-07&1.10 \\
\hline
$\min\{\gamma\sigma,2-\alpha\}$& &1.10&&1.10&&1.10 \\
\end{tabular*}
{\rule{\temptablewidth}{1pt}}
\end{center}\label{table:largeSigma}
\end{table}
     	
\begin{table}[!ht]
\begin{center}
\tabcolsep 0pt {Table 4: \quad Numerical temporal accuracy for $M=600$, $\sigma=0.3$ and $\gamma=5$} \vspace*{0.5pt}
\def\temptablewidth{1.0\textwidth}
{\rule{\temptablewidth}{1pt}}
\begin{tabular*}{\temptablewidth}{@{\extracolsep{\fill}}ccccccc}
 $N$ &\multicolumn{2}{c}{$\alpha=0.3$}&\multicolumn{2}{c}{$\alpha=0.5$}
 &\multicolumn{2}{c}{$\alpha=0.7$} \\
\cline{2-3}    \cline{4-5} \cline{6-7}& $e(M,N)$ & Order($N$) &$e(M,N)$ & Order($N$) &$e(M,N)$&Order($N$)\\
\hline
1024&1.59e-05 &-- &4.60e-05&-- &1.72e-04&-- \\
2048 &5.91e-06 &1.43&1.66e-05&1.47&7.00e-05&1.30\\
4096 &2.11e-06 &1.48&5.95e-06&1.48&2.84e-05&1.30\\
8192&7.43e-07 &1.51&2.13e-06&1.48&1.16e-05&1.30\\
\hline
$\min\{\gamma\sigma,2-\alpha\}$& &1.50&&1.50&&1.30 \\
\end{tabular*}
{\rule{\temptablewidth}{1pt}}
\end{center}\label{table:uniform2}
\end{table}

\begin{table}[!ht]
\begin{center}
\tabcolsep 0pt {Table 5: \quad Numerical temporal accuracy for $M=600$,
$\alpha=0.5$ and $\sigma=0.3$} \vspace*{0.5pt}
\def\temptablewidth{1.0\textwidth}
{\rule{\temptablewidth}{1pt}}
\begin{tabular*}{\temptablewidth}{@{\extracolsep{\fill}}ccccccc}
 $N$ &\multicolumn{2}{c}{$\gamma=4$}
 &\multicolumn{2}{c}{$\gamma=5$}&\multicolumn{2}{c}{$\gamma=6$} \\
\cline{2-3}    \cline{4-5} \cline{6-7}& $e(M,N)$ & Order($N$)
&$e(M,N)$ & Order($N$) &$e(M,N)$&Order($N$) \\
\hline
512 &3.66e-04 &-- &1.26e-04&-- &6.92e-05&-- \\
1024 &1.58e-04 &1.21&4.60e-05&1.46&2.51e-05&1.46\\
2048 &6.89e-05 &1.20&1.66e-05&1.47&9.03e-06&1.47 \\
4096&3.00e-05 &1.20&5.95e-06&1.48&3.24e-06&1.48 \\\hline
$\min\{\gamma\sigma,2-\alpha\}$& &1.20&&1.50&&1.50 \\
\end{tabular*}
{\rule{\temptablewidth}{1pt}}
\end{center}\label{table:largeSigma}
\end{table}

\begin{table}[!ht]
\begin{center}
\tabcolsep 0pt {Table 6: \quad Numerical temporal accuracy for $M=600$,
$\alpha=0.4$ and $\sigma=0.3$ } \vspace*{0.5pt}
\def\temptablewidth{1.0\textwidth}
{\rule{\temptablewidth}{1pt}}
\begin{tabular*}{\temptablewidth}{@{\extracolsep{\fill}}ccccccc}
 $N$ &\multicolumn{2}{c}{$\gamma=4$}
 &\multicolumn{2}{c}{$\gamma=5$}&\multicolumn{2}{c}{$\gamma=6$} \\
\cline{2-3}    \cline{4-5} \cline{6-7}& $e(M,N)$ & Order($N$)
&$e(M,N)$ & Order($N$) &$e(M,N)$&Order($N$) \\	
\hline
512 &2.45e-04 &-- &7.55e-05&-- &3.35e-05&-- \\
1024 &1.10e-04 &1.16&2.75e-05&1.45&1.16e-05&1.53\\
2048 &4.76e-05 &1.21&9.69e-06&1.51&3.99e-06&1.55 \\
4096&2.07e-05 &1.20&3.42e-06&1.50&1.36e-06&1.56 \\\hline
$\min\{\gamma\sigma,2-\alpha\}$& &1.20&&1.50&&1.60 \\
\end{tabular*}
{\rule{\temptablewidth}{1pt}}
\end{center}\label{table:largeSigma}
\end{table}

The numerical result in Table 3 (with $M=600$, $\alpha=0.9$ and $\sigma=1.9$) shows that the discrete scheme of \eqref{scheme-zero-one-5}-\eqref{scheme-zero-one-5} has the temporal order $O(\tau^{2-\alpha})$.
Furthermore, in the case of uniform mesh $\gamma=1$, the solution is accurate of order $O(\tau^{\sigma})$, which matchs with our theoretical analysis of Theorems \ref{the:Convergence-L2}-\ref{the:Convergence-Linfty}.
The numerical results in Tables 4-6 with $M=600$ show that the time accuracy of order $O(\tau^{\min\{\gamma\sigma,2-\alpha\}})$
and support the predicted time accuracy in Theorems \ref{the:Convergence-L2}-\ref{the:Convergence-Linfty}.
The optimal time accuracy $O(\tau^{2-\alpha})$ is observed when the grid parameter $\gamma>(2-\alpha)/\sigma$.
Thus the $L^{2}$ error estimate \eqref{scheme-convegence-zero-one-conclusion-L2-2} and $L^{\infty}$ error estimate \eqref{scheme-convegence-zero-one-conclusion-H2-2} are sharp.




\end{document}